\documentclass[10pt,a4paper]{article}
\usepackage{amsfonts}
\usepackage{amssymb}
\usepackage{color}
\usepackage[centertags, sumlimits, intlimits]{amsmath}
\usepackage[round, sort]{natbib}
\usepackage{float}
\usepackage{cancel}
\usepackage[title]{appendix}

\floatstyle{ruled}
\newfloat{summary}{thp}{lop}
\floatname{summary}{Summary}

\usepackage{listings}

\definecolor{listinggray}{gray}{0.9}
\definecolor{lbcolor}{rgb}{0.9,0.9,0.9}
\definecolor{lbashcolor}{rgb}{0.95,0.99,0.95}
\usepackage{xcolor}
\definecolor{dkgreen}{rgb}{0,0.6,0}
\definecolor{dred}{rgb}{0.545,0,0}
\definecolor{dblue}{rgb}{0,0,0.545}
\definecolor{lgrey}{rgb}{0.9,0.9,0.9}
\definecolor{gray}{rgb}{0.4,0.4,0.4}
\definecolor{darkblue}{rgb}{0.0,0.0,0.6}

\usepackage{listings}
\lstdefinelanguage{cpp}{
    backgroundcolor=\color{lgrey},  
      basicstyle=\scriptsize \ttfamily \color{black} \bfseries,   
      breakatwhitespace=false,       
      breaklines=true,                   
      prebreak = \raisebox{0ex}[0ex][0ex]{\ensuremath{\hookleftarrow}},             
      commentstyle=\color{dkgreen},   
      deletekeywords={...},          
      escapeinside={\%*}{*)},                  
      frame=single,                  
      language=C++,                
      keywordstyle=\color{purple},  
      identifierstyle=\color{black},
      stringstyle=\color{blue},      
      numbers=left,        
      firstnumber=1,                         
      numberstyle=\tiny\color{black}, 
      rulecolor=\color{black},        
      showspaces=false,               
      showstringspaces=false,        
      showtabs=false,                
}

\lstdefinestyle{Bash}
{language=bash,
backgroundcolor=\color{lbashcolor},
keywordstyle=\color{blue},
basicstyle=\scriptsize,
keywordstyle=[2]{\color{red}},
literate={\$}{{\textcolor{red}{\$}}}1 
         {:}{{\textcolor{red}{:}}}1
         {~}{{\textcolor{red}{\textasciitilde}}}1,
}

\usepackage{graphicx,psfrag,epsfig,rotating}
\usepackage{color}

\usepackage{booktabs}
\usepackage{url}
\usepackage{longtable}

\usepackage{multirow}
\usepackage{setspace}
\usepackage{nomencl}
\usepackage{lscape}
\usepackage{tabularx}
\usepackage{pdflscape}

\usepackage{array,float,latexsym,amsmath,amssymb,amsbsy,eurosym,theorem,amsfonts,mathrsfs,bbm,verbatim}

\usepackage[active]{srcltx}

\usepackage{pdfsync}

\newcommand{\fig}[1]{Fig.~\ref{#1}}

\renewcommand{\b}[1]{\boldsymbol{#1}} 

\renewcommand{\c}[1]{\mathcal{#1}}

\renewcommand{\v}[1]{\mathbbm{#1}} 

\newcommand{\h}[1]{\widehat{#1}}
\renewcommand{\o}[1]{\overline{#1}}
\renewcommand{\u}[1]{\underline{#1}}

\newcommand{\lrb}[1]{\left[ {#1} \right]}
\newcommand{\lc}{\left\{}

\newcommand{\p}{\partial}

\newcommand{\f}{\displaystyle\frac}

\renewcommand{\d}{\mbox{d}}
\newcommand{\D}{\mbox{D}}
\newcommand{\grad}{\ensuremath{\mbox{grad}}}

\renewcommand{\div}{\mbox{div}}

\newcommand{\Grad}{\ensuremath{\mbox{Grad}}}

\newcommand{\Div}{\mbox{Div}}

\newcommand{\dyad}{\otimes}

\newcommand{\prs}{{}^{\mbox{\scriptsize p}}} 

\newcommand{\inv}{{}^{-\text{\scriptsize 1}}}
\newcommand{\trns}{{}^{\text{\scriptsize t}}}

\newcommand{\Det}{\text{Det}}

\newcommand{\hsth}{\{\h\bullet\}}

\newcommand{\sth}{\{\bullet\}}

\definecolor{gray}{rgb}{0.75, 0.75, 0.75}
\definecolor{yellow}{rgb}{1, 0.7, 0.2}
\definecolor{green}{rgb}{0.3, 0.9, 0.3}
\definecolor{brown}{rgb}{0.6, 0.3, 0.2}
\definecolor{magenta}{rgb}{0.9, 0.1, 0.9}
\definecolor{light}{rgb}{1, 0.7, 0.7}

\newcommand{\assembly}{\mathop{\mbox{\LARGE {\textsf {\textbf A}}}}}

\title{A finite element implementation of surface elasticity at finite strains using the deal.II library}
\author{Andrew McBride${}^\text{a}$, Ali Javili${}^\text{b}$, Paul Steinmann${}^\text{b}$, Daya Reddy${}^\text{a}$ \\ 
{${}^\text{a}$\footnotesize Centre for Research in Computational and Applied Mechanics, University of Cape Town, South Africa} \\
{${}^\text{b}$\footnotesize Chair of Applied Mechanics, University of Erlangen--Nuremberg, Germany}}
\date{}

\begin{document}
\maketitle

\begin{abstract}
The potentially significant role of the surface of an elastic body in the overall response of the continuum can be described using the mature theory of surface elasticity. 
The objective of this contribution is to detail the finite element approximation of the underlying governing equations (both in the volume and on its surface) and their solution using the open-source finite element library {\tt deal.II}. 
The fully-nonlinear (geometric and material) setting is considered.
The nonlinear problem is solved using a Newton--Raphson procedure wherein the tangent contributions from  the volume and surface are computed exactly.
The finite element formulation is implemented within the total Lagrangian framework and a Bubnov--Galerkin spatial discretization of the volume and the surface employed. 
The surface is assumed material. 
A map between the degrees of freedom on the surface and on the boundary of the volume is used to allocate  the contribution from the surface to the global system matrix and residual vector. 
The  {\tt deal.II} library greatly facilitates the computation of the various surface operators, allowing the numerical implementation to closely match the theory developed in a companion paper.   
Key features of the theory and the numerical implementation are elucidated using a series of benchmark example problems. 
The full, documented source code is provided.
\end{abstract}

\section{Introduction}

The surface elasticity theory of \citet{Gurtin1975} has been widely used to account for the role that the surface of an elastic body can play in the overall response of the continuum. 
Integral to the theory is the derivation of a set of governing equations and constitutive relations that describe the behaviour of the surface of the bulk object.
The role of surface elasticity and the size-dependence of the elastic response has received  considerable attention recently \citep[see e.g.][]{Duan2009, Weissmuller2010}.
This resurgence of interest in the mechanics  of surfaces can be largely attributed to the increasing number of applications involving nanoscale structures.
In such applications scale effects are observed due to the significant surface-to-volume ratio.  
Classical continuum formulations are unable to account for these effects as they lack an inherent length scale.

The theory of surface elasticity is well understood \citep[see e.g.\ the review by][]{Javili2013}. 
The vast majority of numerical treatments of surface elasticity, however, have been limited to the infinitesimal deformation regime.
Furthermore, to the best of the authors knowledge, no open-source or commercial finite element implementation of surface elasticity is available. 

The objective of this contribution is to detail a finite element implementation of surface elasticity at finite deformations using the open-source  library \verb|deal.II| \citep{Bangerth2007, Bangerth, Bangerth2013}. 
The implementation uses various  \verb|deal.II| routines developed to facilitate solving partial differential equations on curved manifolds \citep[see e.g.][]{DeSimone2009, Heltai2008} embedded in a higher-dimensional space.
The  \verb|deal.II| term for this set of routines is \emph{codimension one}. 
The name  follows from the mathematical definition of codimension: if $W$ is a subspace of a linear space or manifold $V$, then the codimension of $W \in V$ is defined by
$\text{codim}\, W = \text{dim}\, V - \text{dim}\, W$. 
The \emph{codimension one} routines greatly simplify the construction of the approximations to the various mathematical  operators on the surface, and the geometric description of the surface. 
This is critical, as the surface deformation gradient (used to parameterize the constitutive response) is rank deficient, and the surface profile can be complex when modelling realistic problems. 
The   \emph{codimension one} routines allow the surface of the body to be treated as an independent two-dimensional manifold embedded in three-dimensional space. 
The constitutive model for the surface implemented here allows the surface to behave in a solid- or fluid-like manner. 
For solid-like behaviour, the surface free energy resembles a classical neo-Hookean model. 
The neo-Hookean model is extended to include surface tension and thereby account for fluid-like effects.

Certain aspects of the implementation presented here for the volume contributions are similar to those discussed in the online tutorial (step\_44) on a three-field formulation for (near-incompressible) finite elasticity (see \url{www.dealii.org/developer/doxygen/deal.II/step_44.html}).

The finite element formulation is implemented within the total Lagrangian framework and a Bubnov--Galerkin spatial discretization of the volume and surface employed. 
The surface is assumed material (i.e.\ it acts like a membrane permanently attached to the underlying solid volume).
The nonlinear problem is solved using a Newton--Raphson procedure wherein the tangent contributions from both the volume and surface are computed exactly. 
A curvilinear-coordinate-based finite element methodology for the problem of surface elasticity was developed in a companion paper \citep{Javili-in-review}.
Extensive theoretical details and references are provided.
For additional information on the admissible range for the surface material parameters, see \citet{Javili2012c}. 
In particular, the validity of negative surface parameters, which have been reported in the literature, is assessed.

While the majority of the numerical examples presented in the literature are performed on academic problems, one intended and novel application of the code presented here is to describe realistic nanostructures with complex surface geometry.
This requires an efficient and parallelized code.  
In order to realise this objective, various numerical operations in the volume and on the surface are parallelized (within a shared memory paradigm) using the Threading Building Blocks library (TBB) \citep{TBB}. 
Aspects of the extension of this framework to a distributed environment is discussed briefly. 

The structure of the paper is as follows. 
The key kinematic concepts are recalled in Section~\ref{sec_kinematics}. 
Thereafter, the independent hyperelastic constitutive relations governing the response of both the volume and the surface are summarised. 
The strong form of the governing equations and their restatement in weak form as a Newton scheme are presented in Section~\ref{sec_governing_equations}. 
The linearized weak form provides the basis for the fully-discrete problem introduced in Section~\ref{sec_fully_discrete}. 
Details of the numerical implementation within the \verb|deal.II| library are given in Section~\ref{sec_numerical_implementation}.
Four numerical examples are presented in Section~\ref{sec_numerical_results}. 
They serve to elucidate key features of the theory and aspects of the numerical implementation.
The complete, documented source code, instructions, and all input decks required to reproduce the numerical examples are provided online at \url{www.cerecam.uct.ac.za/code/surface_energy/doc/html/index.html}.

\subsection*{Notation and definitions}

Direct notation is adopted throughout. 
Occasional use is made of index notation, the summation convention for repeated indices being implied. 
When the repeated indices are lower-case italic letters, the summation is over the range $\{1,2,3\}$.
If they are lower-case Greek letters the summation is over the range $\{1,2\}$.
The scalar product of two vectors $\b{a}$ and $\b{b}$ is denoted $\b{a}\cdot\b{b} = [\b{a}]_{i} [\b{b}]_{i}$.  
The scalar product of two second-order tensors $\b{A}$ and $\b{B}$ is denoted $\b{A}:\b{B} = [\b{A}]_{ij} [\b{B}]_{ij}$. 
The composition  of two second-order tensors $\b{A}$ and $\b{B}$, denoted $\b{A} \cdot \b{B}$, is a second-order tensor with components  $[\b{A} \cdot \b{B}]_{ij} = [\b{A}]_{im} [\b{B}]_{mj}$.
The tensor product of two vectors $\b{a}$ and $\b{b}$ is a second-order tensor $\b{D}=\b{a}\dyad\b{b}$ with $[\b{D}]_{ij}= [\b{a}]_{i} [\b{b}]_{j}$. 
The two non-standard tensor products of two second-order tensors $\b{A}$ and $\b{B}$ are the fourth-order tensors $[ \b{A} \,\o{\dyad}\, \b{B}]_{ijkl} = [\b{A}]_{ik} [\b{B}]_{jl}$ and $[\b{A} \,\u{\dyad}\, \b{B}]_{ijkl} = [\b{A}]_{il} [\b{B}]_{jk}$.

An arbitrary quantity in the volume is denoted as $\sth$ and analogously $\hsth$ denotes an arbitrary surface quantity.
The surface quantity can be a vector, not necessarily tangent to the surface, or a tensor, not necessarily tangential or superficial to the surface.
The (conventional) identity tensor in $\mathsf{E}^3$ is denoted as $\b{i}$ in the spatial configuration and $\b{I}$ in the material configuration.
In what follows the identity tensors $\b{i}$ and $\b{I}$ are understood as the conventional identity tensors in $\mathsf{E}^3$, i.e.\ their matrix representation would be a $3 \times 3$ matrix with 1 in the main diagonal entries and 0 elsewhere. 
Although these identity tensors are invariant and $\b{i} = \b{I}$, we use different letters to indicate explicitly which configuration they belong to.

\section{Summary of the problem of surface elasticity}\label{sec_summary}

The objective of this section is to  recall briefly the problem of surface elasticity at finite deformations. 
For additional details the reader is referred to \citet{Gurtin1975, Javili-in-review, Javili2013}, among others. 

\subsection{Kinematics}\label{sec_kinematics}

\begin{figure}[h]
 \centering
 \includegraphics[width=\textwidth]{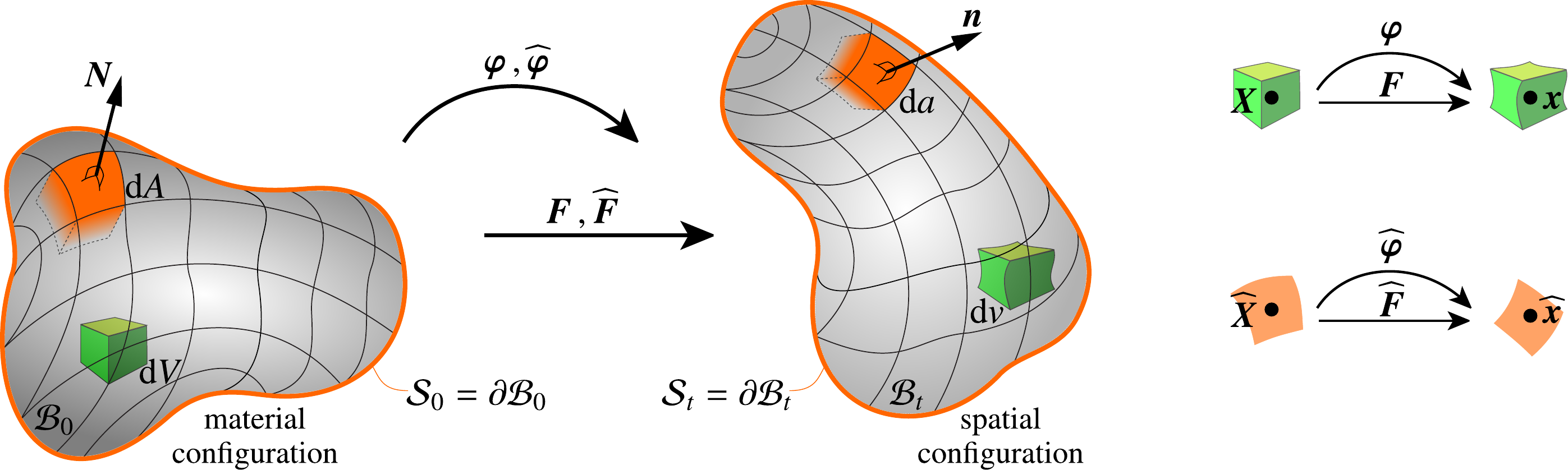}
 \caption{The material and spatial configurations of a continuum body and its boundary and the associated deformation maps and deformation gradients.}
 \label{fig:motion}
\end{figure}

Consider a continuum body that takes the material configuration $\c{B}_0$ at time $t=t_0$ and is mapped via the (volume) non-linear deformation map $\b{\varphi}$ to the spatial configuration $\c{B}_t$ at time $t>0$, as shown in \fig{fig:motion}. 
Material points in the material and spatial configurations are denoted $\b{X}$ and $\b{x}$, respectively.
The associated linear deformation map, i.e.\ the (volume) deformation gradient, is denoted by $\b{F}:=\p \b{x} / \p \b{X} =  {\b{g}}_i \dyad {\b{G}}{}^i$, and maps material line elements $\d\b{X} \in T \c{B}_0$ (tangent to $\c{B}_0$) to spatial line elements $\d{\b{x}}\in T \c{B}_t$ via the relation $\d{\b{x}} = {\b{F}} \cdot \d{\b{X}}$. 
The co- and contravariant basis vectors in the material and spatial configurations are denoted $\b{G}_i$ and $\b{g}^i$, respectively. 
The displacement of a material point is denoted $\b{u} = \b{\varphi}(\b{X},t) - \b{X}$. 
The volume deformation gradient $\b{F}$ is rank-sufficient with inverse $\b{f}:=\p \b{X} / \p \b{x}= {\b{G}}_i \dyad {\b{g}}{}^i$. 
The determinant of the deformation gradient and its inverse are given by $J = \det \b{F} > 0$ and  $j = \det \b{f} = 1 / J > 0$, respectively.

Let $\c{S}_0$ and $\c{S}_t$ denote the surface of the continuum body in the material and spatial configurations, respectively.
The outward unit normals to $\c{S}_0$ and $\c{S}_t$ are denoted $\b{N}$ and $\b{n}$, respectively. 
Material particles on the surface are denoted $\h{\b{X}}$ in the material configuration and are \emph{attached} to the volume, i.e.\ $\h{\b{X}}=\b{X}|_{\p{\c{B}_0}}$.
Consequently, $\c{S}_0 = \p\c{B}_0$.
Furthermore, we assume that the surface is material in the sense that it is permanently attached to the substrate (i.e.\ the boundary of the volume).
Therefore $\c{S}_t = \p\c{B}_t$ and $\h{\b{x}}=\b{x}|_{\p{\c{B}_t}}$.
This assumption implies that the motion of the surface $\h{\b{\varphi}}$ is the restriction of the volume motion $\b{\varphi}$ to the surface, i.e.\ $\h{\b{\varphi}} = \b{\varphi}|_{\p{\c{B}_0}}$.

Material line elements on the surface in the material and spatial configurations are denoted  $\d\h{\b{X}} \in T \c{S}_0$ and $\d\h{\b{x}} \in T \c{S}_t$, respectively, and are related as $\d\h{\b{x}} = \h{\b{F}} \cdot \d\h{\b{X}}$, where $\h{\b{F}}:=\p \h{\b{x}} / \p \h{\b{X}}  = {\h{\b{g}}}_\alpha \dyad {\h{\b{G}}}{}^\alpha$ denotes the rank-deficient, and thus non-invertible, surface deformation gradient.
The co- and contravariant surface basis vectors in the material and spatial configurations are denoted $\h{\b{G}}_\alpha$ and $\h{\b{g}}^\alpha$, respectively. 
Various key concepts from differential geometry are given in Summary~\ref{summary:surface-1}.
The inverse deformation gradient $\h{\b{f}}$ is defined by $\h{\b{f}}:=\p \h{\b{X}} / \p \h{\b{x}}  = {\h{\b{G}}}_\alpha \dyad {\h{\b{g}}}{}^\alpha$ and is related to the surface deformation gradient as follows:
\begin{align*}
    \h{\b{f}} \cdot \h{\b{F}} = \h{\b{I}} =: \b{I} - \b{N}\otimes\b{N}
    && \text{and} &&
    \h{\b{F}} \cdot \h{\b{f}} = \h{\b{i}} =: \b{i} - \b{n}\otimes\b{n} \, .
\end{align*}
The (surface) determinant of the surface deformation gradient and its inverse are denoted $\h{J} = \h{\det} \h{\b{F}} > 0$ and $\h{j} = \h{\det} \h{\b{f}} > 0$, respectively, where $\h{j} = 1 / \h{J}$. 

\begin{summary}[h]

\begin{minipage}{0.3\textwidth}
 \centering
 \includegraphics[width=1.0\textwidth]{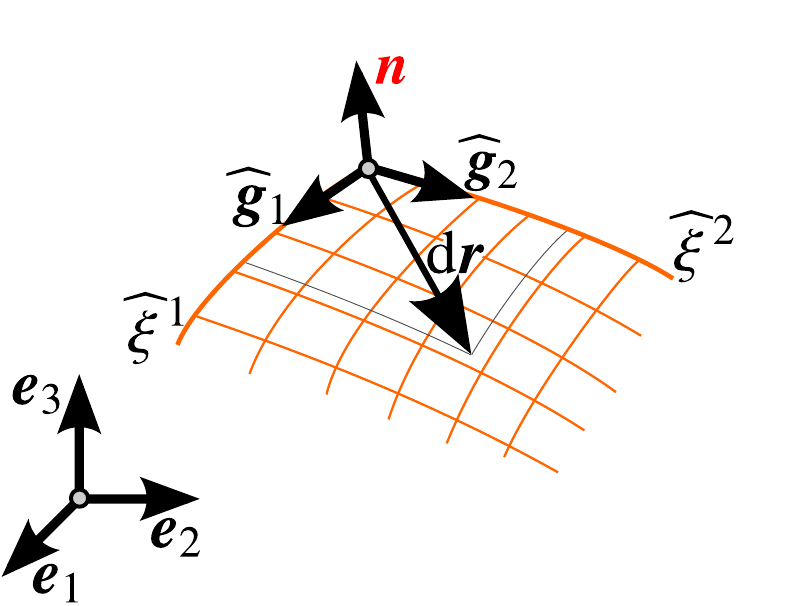}
\end{minipage}
\hfill
\begin{minipage}{0.65\textwidth}
  \begin{equation*}
    \begin{array}{c}
    \d\b{r} = \d\b{r} (\h{\b{\xi}}) = \d\b{r} (\h{\xi}{}^1,\h{\xi}{}^2)\\[12pt]
    \h{\b{g}}_\alpha = \f{\p\b{r}}{\p\h{\xi}^\alpha} \qquad , \qquad \h{\b{g}}^\alpha = \f{\p\h{\xi}^\alpha}{\p\b{r}} \qquad \text{with} \qquad \alpha \in \{ 1,2 \} \\[15pt]
    \h{\b{g}}{}^3 = \h{\b{g}}_1 \times \h{\b{g}}_2 \quad , \quad \h{\b{g}}_3 = [\h{g}{}^{33}]\inv \, \h{\b{g}}{}^3 \quad , \quad \b{n} = \sqrt{\h{g}{}_{33}} \, \h{\b{g}}^3 = \sqrt{ \h{g}{\,}^{33}} \, \h{\b{g}_3} \\[15pt]
    \h{\grad}\hsth = \f{\p\hsth}{\p\h{\xi}{}^\alpha} \dyad \h{\b{g}}_\alpha \qquad , \qquad \h{\div}\hsth = \f{\p\hsth}{\p\h{\xi}{}^\alpha} \cdot \h{\b{g}}_\alpha = \h{\grad}\hsth : \h{\b{i}} \\[15pt]
    \h{\det}\hsth = \f{\big|[\hsth\cdot\h{\b{g}}_1]\times[\hsth\cdot\h{\b{g}}_2]\big|}{| \h{\b{g}}_1 \times \h{\b{g}}_2|} 
 \end{array}
  \end{equation*}
\end{minipage}

  \begin{equation*}
    \begin{array}{c}
    \h{\b{g}}_\alpha = \h{g}_{\alpha\beta} \, \h{\b{g}}^\beta \, , \quad \h{g}_{\alpha\beta} = \h{\b{g}}_\alpha \cdot \h{\b{g}}_\beta \, , \qquad \h{\b{g}}{}^\alpha = \h{g}\,{}^{\alpha\beta} \, \h{\b{g}}_\beta \, , \quad \h{g}\,{}^{\alpha\beta} = \h{\b{g}}^\alpha \cdot \h{\b{g}}^\beta \, , \qquad [\h{g}_{\alpha\beta}]=[\h{g}\,{}^{\alpha\beta}]\inv \, , \qquad \h{g} = |[\h{g}_{\alpha\beta}]|  \\[10pt]
    \h{\v{g}}: \text{surface permutation tensor} \, , \, \h{\v{g}} = \h{g}_{\alpha\beta} \, \h{\b{g}}{}^\alpha\dyad\h{\b{g}}{}^\beta\dyad\b{g}^3 = \h{g}\,{}^{\alpha\beta} \, \h{\b{g}}_\alpha\dyad\h{\b{g}}_\beta\dyad\b{g}_3   \, , \qquad \h{g}_{\alpha\beta} = \h{g}\,{}^{\alpha\beta} \ \h{g} = \sqrt{\, \h{g}}\; \h{e}_{\alpha\beta} \\[12pt]
    \h{e}_{\alpha\beta} = \lc \begin{array}{cl}
                           1 & \text{if $\alpha\beta$ is 12}\\
                           -1 & \text{if $\alpha\beta$ is 21}\\
                           0 & \text{otherwise}
                          \end{array} \right. \, , \quad  \h{g}_{\alpha\beta} = | [\h{\b{g}}_\alpha \times \h{\b{g}}_\beta]| = \sqrt{\, \h{g}}\; \h{e}_{\alpha\beta} \, , \quad  \h{g}\,{}^{\alpha\beta} = |[\h{\b{g}}{}^\alpha \times \h{\b{g}}{}^\beta]| = [\sqrt{\, \h{g}}]\inv\; \h{e}_{\alpha\beta}\\[18pt]
    \h{\b{u}} \times \h{\b{v}} = [ \h{\b{u}} \dyad \h{\b{v}}] : \h{\v{g}} \, , \, \h{\b{u}} \cdot \h{\b{v}} = [ \h{\b{u}} \dyad \h{\b{v}}] : \h{\b{i}} \, , \,
    \h{\b{i}} = {\delta}^\alpha_\beta \, \h{\b{g}}_\alpha\dyad\h{\b{g}}{}^\beta = \h{\b{g}}_\alpha\dyad\h{\b{g}}{}^\alpha = \h{\b{g}}_1\dyad\h{\b{g}}{}^1 + \h{\b{g}}_2\dyad\h{\b{g}}{}^2 = \b{i} - \b{n}\dyad\b{n} 
 \end{array}
  \end{equation*}

  \caption{The key differential geometry concepts of the surface \citep[see][]{Javili-in-review}. The surface coordinates are denoted $\h{\b{\xi}}$.}
  \label{summary:surface-1}
\end{summary}

\subsection{Constitutive relations}

A hyperelastic, neo-Hookean material model is assumed for the constitutive response of the volume (see Summary~\ref{constiutive_relations}). 
The free energy $\Psi = \Psi(\b{F})$ is parametrised by the Lam\'{e} moduli $\lambda$ and $\mu$. 
The constitutive relation on the surface is chosen to mimic that in the volume. 
In addition to the neo-Hookean type surface hyperelastic response, the surface free energy $\h{\Psi} = \h{\Psi}(\h{\b{F}})$ accounts for surface tension via the parameter $\h{\gamma}$.  
The contribution of surface tension renders the surface free energy non-zero in the material configuration. 
This has implications for the resulting numerical scheme.
For a study of the admissible range for the surface material parameters, see \citet{Javili2012c}.

\begin{summary}[htb]
\textbf{Volume}:
\begin{align*}
 J := \Det{\b{F}} \quad \text{: Jacobian determinant} && \mu \,, \lambda \quad  \text{: Lam\'{e} constants} 
\end{align*}
 \begin{align*}
  &\text{Free energy} && \Psi(\b{F}) = \tfrac{1}{2} \, \lambda \ln^2 J + \tfrac{1}{2} \, \mu \, [\b{F}:\b{F} - 3 - 2 \ln J] \\[8pt]
  &\text{Piola--Kirchhoff stress} && \b{P}(\b{F}) = \f{\p\Psi}{\p\b{F}} = \lambda \ln J \, \b{f}\trns + \mu \, [\b{F} - \b{f}\trns] \\[8pt]
  &\text{Piola stress tangent} && \v{A}(\b{F}) = \f{\p\b{P}}{\p\b{F}} = \lambda\, [\b{f}\trns \dyad \b{f}\trns + \ln J \, \v{D}] + \mu \, [\v{I} - \v{D}]
 \end{align*}
\begin{align*}
 \v{D} := \f{\p\b{f}\trns}{\p\b{F}} = - \b{f}\trns \, \u\dyad \, \b{f} \qquad , \qquad \v{I} := \f{\p\b{F}}{\p\b{F}} = \b{i} \, \o\dyad \, \b{I}
\end{align*}
\hrulefill \\
\textbf{Surface}:
\begin{align*}
 \h{J} := \h{\Det}{\h{\b{F}}} \quad \text{: Surface Jacobian determinant} && \h{\mu} \,, \h{\lambda} \quad  \text{: Surface Lam\'{e} constants} && \h{\gamma} \quad  \text{: Surface tension}  
\end{align*}
 \begin{align*}
  &\text{Surface Free energy} && \h{\Psi}(\h{\b{F}}) = \tfrac{1}{2} \, \h{\lambda} \ln^2 \h{J} + \tfrac{1}{2} \, \h{\mu} \, [\h{\b{F}}:\h{\b{F}} - 2 - 2 \ln \h{J}] + \h{\gamma} \, \h{J} \\[8pt]
  &\text{Surface Piola--Kirchhoff stress} && \h{\b{P}}(\h{\b{F}}) = \f{\p\h{\Psi}}{\p\h{\b{F}}} = \h{\lambda} \ln \h{J} \ \h{\b{f}}\trns + \h{\mu} \, [\h{\b{F}} - \h{\b{f}}\trns] + \h{\gamma}\, \h{J}\ \h{\b{f}}\trns \\[6pt]
  &\text{Surface Piola stress tangent} && \h{\v{A}}(\h{\b{F}}) = \f{\p\h{\b{P}}}{\p\h{\b{F}}} = \h{\lambda}\, [\h{\b{f}}\trns \dyad \h{\b{f}}\trns + \ln \h{J} \, \h{\v{D}}] + \h{\mu} \, [\h{\v{I}} - \h{\v{D}}] + \h{\gamma}\, \h{J}\ [ \h{\b{f}}\trns \dyad \h{\b{f}}\trns + \h{\v{D}}]
 \end{align*}
\begin{align*}
 \h{\v{D}} := \f{\p\h{\b{f}}\trns}{\p\h{\b{F}}} = - \h{\b{f}}\trns \, \u\dyad \, \h{\b{f}} + [\b{i}-\h{\b{i}}]\, \o\dyad \, [\h{\b{f}}\cdot\h{\b{f}}\trns] \qquad , \qquad \h{\v{I}} := \f{\p\h{\b{F}}}{\p\h{\b{F}}} = \b{i} \, \o\dyad \ \h{\b{I}}
\end{align*}
 \caption{Hyperelastic material models for the volume and surface}
 \label{constiutive_relations}
\end{summary}

\subsection{Governing equations}\label{sec_governing_equations}

The strong form of the equations governing the response of the volume and the surface (i.e.\ the quasi-static balances of linear momentum) are given by
\begin{align}
\Div \b{P} + \b{b}_0^\text{p} = \b{0} && \text{in } \c{B}_0 \, , \label{strong_bulk} \\
\h{\Div} \h{\b{P}} + \h{\b{b}}_0^\text{p} - \b{P}\cdot\b{N} = \b{0} && \text{on } \c{S}_0 \, , \label{strong_surface}
\end{align}
where $\b{b}_0^\text{p}$ is the (prescribed) body force per unit reference volume in the material configuration, and $\h{\b{b}}_0^\text{p}$ is the (prescribed) force per unit reference area of the material configuration. 
Note that in the absence of a surface with an independent free energy (i.e.\ an \emph{energetic surface}), (\ref{strong_surface}) defines the standard Neumann boundary condition on the Piola traction, i.e.\ $\h{\b{b}}_0^\text{p} \equiv \b{P}\cdot\b{N}$. 

The weak form is obtained by respectively testing (\ref{strong_bulk}) and (\ref{strong_surface}) with arbitrary test functions $\delta \b{\varphi} \in H^1_0(\c{B}_0)^3$ and $\delta \h{\b{\varphi}} \in H^1_0(\c{S}_0)^3$, applying the divergence theorem in its extended form \citep[see e.g.][]{Javili2013}, exploiting the orthogonality properties of the surface Piola--Kirchhoff stress measure (i.e.\ $\h{\b{P}}\cdot \b{N} = \b{0}$), and the kinematic constraint on the motion of the surface $\delta \h{\b{\varphi}} = \delta \b{\varphi}|_{\partial \c{B}_0}$, rendering the weak form as:
\begin{equation}\label{weak_combined}
\begin{split}
 0 &= \int_{\c{B}_0} \Grad  \delta\b{\varphi} : \b{P}(\b{F})  \, \d V 
 	 + \int_{\c{S}_0} \h{\Grad} \delta \h{\b{\varphi}} :\h{\b{P}}(\h{\b{F}})  \, \d A
	- \int_{\c{B}_0}  \delta \b{\varphi} \cdot \b{b}\prs_0 \, \d V 
	- \int_{\c{S}_0^\text{\tiny N}} \delta \h{\b{\varphi}} \cdot \h{\b{b}}\prs_0 \, \d A \\
	&\quad-\cancelto{0}{\biggl[\int_{\partial \c{B}_0} \delta\b{\varphi} \cdot [\b{P} \cdot \b{N}] \, \d A - \int_{\c{S}_0} \delta\h{\b{\varphi}} \cdot [\b{P} \cdot \b{N}] \, \d A \biggr]}
	\, , 
\end{split}	
\end{equation}
where the dependence of the Piola--Kirchhoff stress measures $\b{P}$ and $\h{\b{P}}$  on the solution is given via the constitutive relations in Summary~\ref{constiutive_relations}. 
The Neumann part of the surface is denoted $\c{S}_0^\text{\tiny N}$. 
Note that from the assumption of a material surface $\delta\b{\varphi}\vert_{\c{S}_0} = \delta\h{\b{\varphi}}$, the final term of (\ref{weak_combined}) is thus zero due to traction continuity (as indicated).
It will prove convenient, when constructing the finite element approximation, to decompose (\ref{weak_combined}) into  volume and surface contributions. 
The volume and surface contributions will be recombined when the linear system is assembled. 

The decomposed weak form of the governing equations is solved using a Newton--Raphson strategy. 
The time domain $[0,T]$ is decomposed into $N$ uniform intervals of duration $\Delta t:= T / N = t_{n+1} - t_n$, where $t_n = n \Delta t$.
The complete state of the system is assumed known at $t_n$.
The (not necessarily converged) value of a variable evaluated at an iteration $(i)$ during the $n+1$ timestep is denoted $(\bullet)_{n+1}^{(i)} \equiv (\bullet)^{(i)}$.
The volume and surface residual contributions, denoted $R$ and $\h{R}$, respectively, and the directional derivatives of the residual contributions in the direction of the solution increment in the volume  $\Delta \b{u} = \b{\varphi}^{(i+1)} - \b{\varphi}^{(i)}$  and on the surface $\Delta \h{\b{u}}$ are defined by
\begin{align*}
 R^{(i)} &:= \int_{\c{B}_0} \Grad  \delta\b{\varphi} : \b{P}^{(i)}  \, \d V 
 	- \int_{\c{B}_0}  \delta \b{\varphi} \cdot  \b{b}\prs_0 \, \d V \, , \\
\h{R}^{(i)}&:=	\int_{\c{S}_0} \h{\Grad} \delta \h{\b{\varphi}} :\h{\b{P}}{}^{(i)}  \, \d A 
	- \int_{\c{S}_0^\text{\tiny N}} \delta \h{\b{\varphi}} \cdot \h{\b{b}}\prs_0 \, \d A \, ,
\intertext{and}
 \D_{\Delta \b{u}} R^{(i)} &:= \int_{\c{B}_0} \Grad  \delta\b{\varphi} : \dfrac{\p \b{P}{}^{(i)}}{\p \b{F}} : \Grad{\Delta \b{u}}  \, \d V \, , \\
\D_{\Delta \h{\b{u}}}\h{R}^{(i)} &:=\int_{\c{S}_0} \h{\Grad} \delta \h{\b{\varphi}} : \dfrac{\p \h{\b{P}}{}^{(i)}}{\p \h{\b{F}}} : \h{\Grad} \Delta \h{\b{u}}   \, \d A  \, . 
\end{align*}
The derivatives of the Piola--Kirchhoff stress measures with respect to the deformation measures (i.e.\ the tangents) are given in Summary~\ref{constiutive_relations}. 
Note, an exact computation of the tangent in the volume and on the surface is used.

The resulting Newton scheme is thus given by
\begin{align}
 R^{(i)} + \D_{\Delta \b{u}}R^{(i)}  + \h{R}^{}{(i)} + \D_{\Delta \h{\b{u}}}\h{R}{}^{(i)} & \equiv 0 \, , \notag \\
\implies \qquad \D_{\Delta \b{u}}R^{(i)}   + \D_{\Delta \h{\b{u}}}\h{R}{}^{(i)} & =  -R^{(i)} - \h{R}^{}{(i)} \, .
  \label{linear_combined}
\end{align}

\section{The fully-discrete problem}\label{sec_fully_discrete}

\begin{figure}[!ht]
\centering
\includegraphics[width = 0.9\textwidth]{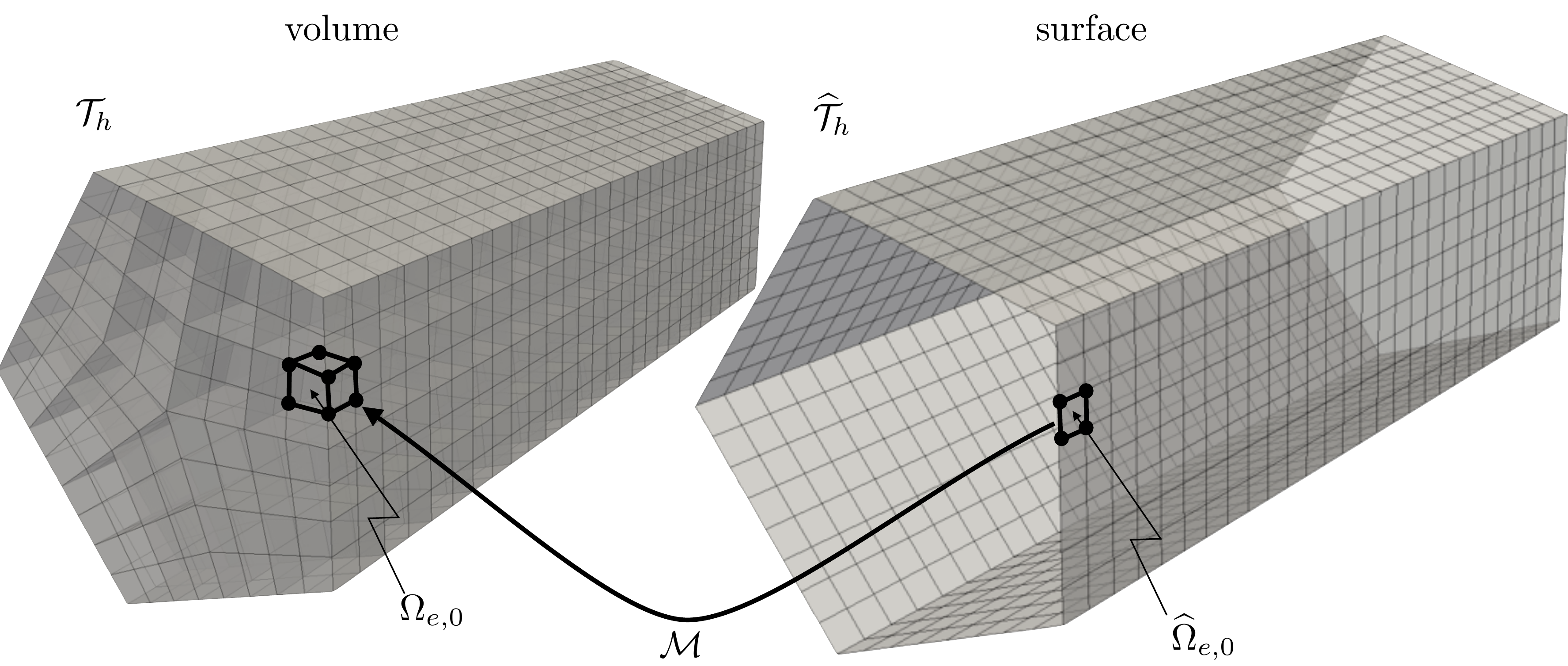}
 \caption{The triangulation of the material volume $\mathcal{T}_h$ and a typical volume cell $\Omega_{e,0}$. 
 The surface triangulation $\widehat{\mathcal{T}}_h$ is extracted from the volume triangulation. 
 A typical cell on the surface is denoted $\widehat{\Omega}_{e,0}$.}
 \label{triangulation}
\end{figure}

The material volume and the surface are partitioned into sets of non-overlapping cells (elements) individually denoted $\Omega_{e,0}$ and $\h{\Omega}_{e,0}$, respectively (see Fig.~\ref{triangulation}).  
The discrete form of the governing equations is obtained by approximating the test functions and trial solutions in the linearized weak form (\ref{linear_combined}) using vector-valued shape functions.
The volume and surface shape functions, associated with an arbitrary degree of freedom in the volume $I$ and on the surface $\h{I}$, are denoted $\b{\Phi}^I$ and $\h{\b{\Phi}}{}^{\h{I}}$, respectively. 
The associated value of the degree of freedom (i.e.\ the increment in the motion) in the volume and on the surface are denoted $\mathsf{u}^I$ and $\h{\mathsf{u}}{}^{\h{I}}$, respectively. 
They form the entries of the global solution vectors $\mathsf{u}$ and $\h{\mathsf{u}}$. 
The trial solution in the volume and on the surface, and their gradients are approximated as follows\footnote{Here and henceforth we adopt the abbreviated summation notation: $\sum_{I} \equiv \sum_{I}^{\text{n}_{\text{dof}}}$ where $\text{n}_{\text{dof}}$ is the number of degrees of freedom in the volume (or on the surface).}
\begin{align*}
\b{u} \approx \sum_I \b{\Phi}^I \mathsf{u}^I	&& \text{and}
	&& \Grad \b{u} \approx \sum_I \Grad \b{\Phi}^I \mathsf{u}^I \, , \\
\h{\b{u}} \approx \sum_{\h{I}} \h{\b{\Phi}}{}^{\h{I}} \h{\mathsf{u}}^{\h{I}}	&& \text{and}
	&& \h{\Grad} \h{\b{u}} \approx \sum_{\h{I}} \h{\Grad} \h{\b{\Phi}}{}^{\h{I}} \h{\mathsf{u}}^{\h{I}} \, .
\end{align*} 
The same approximations are used for the test functions (i.e.\ a Bubnov--Galerkin spatial discretization is employed).

The resulting linearised incremental problem to be solved is given in matrix form as 
\begin{align}
\underbrace{\begin{bmatrix}
\mathsf{A}^{(i)} + \h{\mathsf{A}}^{(i)}
\end{bmatrix}}_{\overset{\square}{\mathsf{A}}}
\begin{bmatrix}
\Delta \mathsf{u}
\end{bmatrix}
=
-
\underbrace{\begin{bmatrix}
 \mathsf{R}^{(i)} + \h{\mathsf{R}}^{(i)}
\end{bmatrix}}_{\overset{\square}{\mathsf{R}}} \, , \label{newton_step}
\end{align}
where the global residual vector $\overset{\square}{\mathsf{R}}$ and system matrix $\overset{\square}{\mathsf{A}}$ are assembled from the element contributions as follows:
\begin{align*}
\mathsf{R}^{(i)} = \assembly_e \sum_{I} \mathsf{R}^I_e
&& \text{and} &&
\mathsf{A}^{(i)} = \assembly_e \sum_{I}\sum_{J} \mathsf{A}^{IJ}_e \, , \\
\h{\mathsf{R}}^{(i)} = \assembly_e \sum_{\h{I}} \h{\mathsf{R}}^{\h{I}}_e
&& \text{and} &&
\h{\mathsf{A}}^{(i)} = \assembly_e \sum_{\h{I}}\sum_{\h{J}} \h{\mathsf{A}}^{\h{I}\h{J}}_e \, ,
\end{align*}
where
\begin{align*}
 \mathsf{R}^I_e &:= \int_{\Omega_{e,0}} \Grad  \b{\Phi}^I : \b{P}  \, \d V 
 	- \int_{\Omega_{e,0}}   \b{\Phi}^I \cdot \b{b}\prs_0 \, \d V \, , \\
 \h{\mathsf{R}}^{\h{I}}_e &:=	\int_{\h{\Omega}_{e,0}} \h{\Grad} \h{\b{\Phi}}{}^{\h{I}} :\h{\b{P}}  \, \d A 
	- \int_{\Omega_{e,0}^\text{\tiny N}}  \h{\b{\Phi}}{}^{\h{I}} \cdot \h{\b{b}}\prs_0 \, \d A \, , \\
 \mathsf{A}^{IJ}_e &:= \int_{\Omega_{e,0}} \Grad  {\b{\Phi}}^I : \v{A} : \Grad {\b{\Phi}}^J \, \d V \, , \\
  \h{\mathsf{A}}^{\h{I}\h{J}}_e &:= \int_{\Omega_{e,0}} \h{\Grad}  {\h{\b{\Phi}}}{}^{\h{I}} : \h{\v{A}} : \h{\Grad} \h{\b{\Phi}}{}^{\h{J}} \, \d A \, ,
\end{align*}
and $\displaystyle{\assembly_e}$ denotes the (non-standard) assembly operator. 
The assembly operator is defined such that contributions from degrees of freedom on the surface are mapped to the corresponding degree of freedom in the volume (the surface is material and the triangulation of the surfaces matches the triangulation of the boundary of the volume). 
The number of surface degrees of freedom will always be less than those in the volume.
The assembly operator is defined such that the size of the matrices $\mathsf{A}^{(i)}$ and $\h{\mathsf{A}}^{(i)}$ are the same.
In other words $\h{\mathsf{A}}^{(i)}$ will contain empty rows and columns corresponding to degrees of freedom internal to the volume.

\section{Numerical implementation using the finite element library {\tt deal.II}}\label{sec_numerical_implementation}

The open-source finite element library {\tt deal.II} \citep{Bangerth2007, Bangerth} is used to assemble and solve the discrete system of equations governing the problem of surface elasticity.  
The implementation here uses various {\tt deal.II} routines developed for the solution of partial differential equations on curved manifolds \citep[see e.g.][]{DeSimone2009, Heltai2008}.

The objective of this section is to review various key features of the numerical implementation, and to discuss some of the overarching design decisions. 
The complete, documented source-code and instructions can be found online at  \url{www.cerecam.uct.ac.za/code/surface_energy/doc/html/index.html}. 
Certain aspects of the implementation are similar to those discussed in the online tutorial (step\_44) on (near-incompressible) finite elasticity without surface effects (see \url{www.dealii.org/developer/doxygen/deal.II/step_44.html}).

\subsection{Triangulations and degree of freedom handlers}

The triangulation (a class in {\tt deal.II}) of the volume, denoted $\mathcal{T}_h$, consists of the location of the vertices of the volume cells (elements) and the cell connectivity, as depicted in Fig.~\ref{triangulation}.
A degree of freedom handler (another class in {\tt deal.II}) combines the purely geometrical information of the triangulation with details of the finite element interpolation space. 
The surface degree of freedom handler $\widehat{\mathcal{D}}_h$ is extracted directly from the volume mesh $\mathcal{D}_h$ using the function \verb|GridTools::extract_boundary_mesh|.
The surface degree of freedom handler is independent to that of the volume. 
A map $\mathcal{M}$, named \verb|surface_to_volume_dof_map| in the code, is used to link the degrees of freedom on the surface to those in the volume.
The surface triangulation is denoted $\widehat{\mathcal{T}}_h$.

The contribution from the volume to the global system matrix is obtained by looping over all the volume cells $\Omega_{e,0}$, assembling the matrix contribution from the cell $\mathsf{A}_e$, and then adding this to the global system matrix. 
This is governed by the method {\tt Solid<...>::assemble\_system\_volume}.
The contribution from the surface is obtained in a similar way.
A loop is performed over all surface cells $\h{\Omega}_{e,0}$ and the tangent contribution  $\h{\mathsf{A}}_e$ calculated. 
The map $\mathcal{M}$ between the degrees of freedom on the surface and the corresponding ones in the volume  is then used to add the surface cell contributions  to the global system matrix (recall the surface is material). 

The Newton scheme (\ref{newton_step}) continues until the normalised measure (relative to the first iteration of the current time step) of the solution increment $\Delta \mathsf{u}$ and the residual $\overset{\square}{\mathsf{R}}$ decreases below a user-specified tolerance. 

It is important to note that alternatively a separate global surface system matrix could be constructed as the degrees of freedom on the surface are independent to those of the volume. 
The degrees of freedom on the surface could then be constrained to be the same as those in the volume and the resulting global block system solved. 
Using this approach, one could also then consider the case of non-material surfaces. 
It was decided that this approach would be more cumbersome as the global system matrix would have to be initialized using sparsity patterns from different triangulations. 
The approach adopted here can still be used to solve membrane problems (i.e.\ where the volume is absent) by defining a volume with vanishing strength.
More details are given in the liquid bridge numerical example  discussed in Sect.~\ref{sec:liquid_bridge}. 
It should be emphasised, however, that the solution of membrane problems is not the focus of this contribution.

\subsection{The {\tt ContinuumPoint} class}\label{sec:continuum_points}

The \verb|ContinuumPoint| class represents a Lagrangian continuum point in the volume or on the surface. 
Its primary role is to calculate the constitutive response at a quadrature point given the kinematic state. 
For example, it determines the Piola--Kirchhoff stress and the tangent based on the deformation gradient, and its determinant and inverse. 
It is also convenient to store the kinetic and kinematic data in the class structure for use when constructing the residual and for post-processing the results. 
One key design decision was not to create separate classes for the material response and the quadrature point data as done in step\_44 as this leads to data ownership conflicts.

An instance of a \verb|ContinuumPoint|  is aware if it is in the volume or on the surface. 
Such knowledge is required occasionally (for example, when computing the tangent) but for the majority of the operations it is not required (for example, when computing the stress).

\subsection{Computation of kinematic quantities}

\begin{figure}[!ht]
\centering
\includegraphics[width = \textwidth]{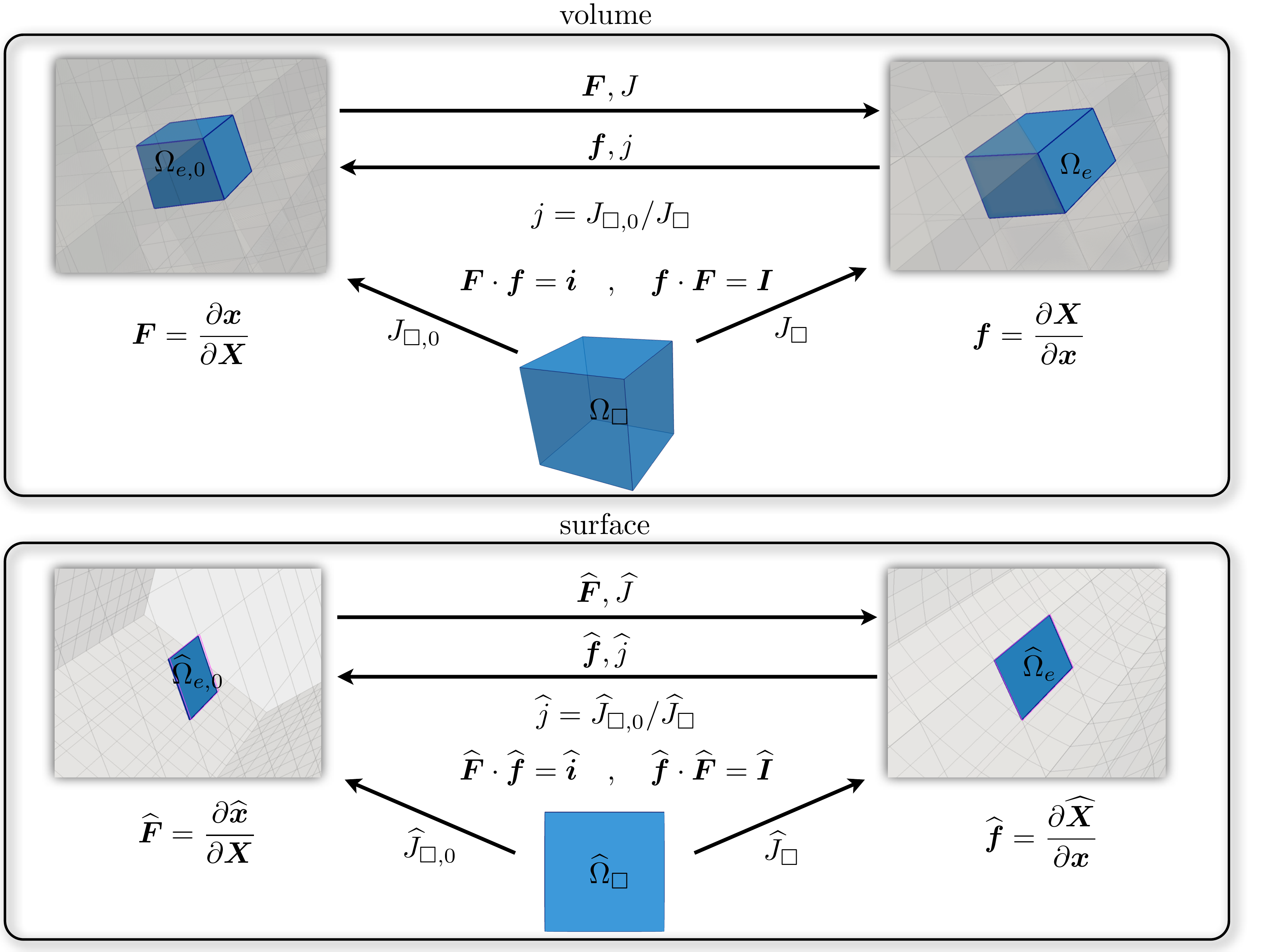}
 \caption{The kinematics of the motion in the volume and on the surface. 
 Also shown in the Lagrangian and Eulerian view of a cell relative to the reference cell.}
 \label{kinematics}
\end{figure}

The form of the constitutive response at a quadrature point within the volume is near identical to that on the surface (see Summary \ref{constiutive_relations}). 
The primary difference is  the structure of the deformation gradient (and its inverse) used to parametrise the free energy. 
The computation of the deformation gradient in the volume and on the surface are performed here in an identical manner (as depicted in Fig.~\ref{kinematics}) using the routines \verb|Solid<spacedim>::update_volume_cp_incremental_one_cell| and \verb|Solid<spacedim>::update_surface_cp_incremental_one_cell|.

The reference volume cell (i.e.\ the isoparametric domain) is denoted $\Omega_\square$. 
A typical volume cell in the material configuration is denoted  $\Omega_{e,0}$. 
A spatial view of the volume cell, denoted $\Omega_{e}$,  is obtained by applying the nonlinear motion map $\b{\varphi}$ to all points with coordinates $\b{X} \in \Omega_{e,0}$. 
This is done in the code using the class \verb|MappingQEulerian|.
The material and spatial gradients of an arbitrary field can then be interpolated from the material and spatial views of the cell as follows:
\begin{align*}
	\Grad (\bullet) \approx \sum_{I} [\Grad{\b{\Phi}^I(\b{X})}][(\bullet)^I]  
	&& \text{and} &&
	\grad (\bullet) \approx \sum_{I} [\grad{\b{\Phi}_t^I(\b{x})}][(\bullet)^I] \, ,
\end{align*}
where a shape function in the spatial configuration associated with an arbitrary node $I$ is denoted $\b{\Phi}^I_t$.
By choosing the arbitrary field in the material view as the spatial placement $\b{x}$ one obtains $\b{F}$, while
choosing the arbitrary field in the spatial view as the material placement $\b{X}$ one obtains $\b{f}$. 

The volume Jacobian determinant $J$ and its inverse $j$ are determined from the ratio of the Jacobian determinant of the map from the isoparametric configuration to the reference configuration $J_{\square,0}$ to the inverse of the Jacobian determinant of the map from the isoparametric configuration (reference cell) to the spatial configuration $J_{\square}$.

An identical approach to the volume is followed on the surface. 
The functionality provided by the \emph{codimension routines} in {\tt  deal.II} allows the surface gradient of a field to be evaluated in a straightforward manner. 
The surface deformation gradient and its inverse are thus computed without needing to invert $\h{\b{F}}$. 
This is critical as $\h{\b{F}}$ is rank deficient. 
The surface Jacobian determinant $\h{J}$ and its inverse $\h{j}$ are determined from the ratio of the Jacobian determinant of the map from the isoparametric configuration to the reference configuration $\h{J}_{\square,0}$ to the inverse of the Jacobian determinant of the map from the isoparametric configuration to the spatial configuration $\h{J}_{\square}$.

\subsection{Parallelization of key tasks}

The  library TBB \citep{TBB} is used to perform as many computationally intensive distributed tasks as possible. 
These include:
\begin{enumerate}
\item the assembly of the tangent and residual contributions from the cells in the volume and on the surface;
\item updating the state of the \verb|ContinuumPoint| after solving the linearized problem in the Newton scheme.
\end{enumerate}
The main tool for using TBB in {\tt deal.II} is the {\tt WorkStream} class.
The bottleneck, however, for large computations is the linear solver. 
The code could, in theory, be parallelized in a distributed sense relatively easily. 
The primary challenge would be to ensure that the independent distributed volume and surface meshes can communicate, via a map of the common degrees of freedom.

\subsection{The choice of linear solver and preconditioner}

The choice of linear solver and preconditioner are specified in the parameter file \verb|parameter.prm|. 
The default choice is the conjugate gradient solver (the matrix problem is symmetric) with Jacobi preconditioning as provided by the {\tt deal.II} library. 
Both this choice of preconditioner and solver are multithreaded. 
The choice works well for the majority of problems investigated.

\section{Numerical results}\label{sec_numerical_results}

The objective of the current section is to elucidate key features of the formulation using four example problems. 
The first example is the Cook's membrane problem. 
This is a widely-used benchmark. 
The purpose of the example is to provide detailed information for a series of uniformly refined meshes. 
Neo-Hookean type surface effects are considered. 
The second example illustrates neo-Hookean type surface effects in a nanowire undergoing significant tensile extension.  
Surface tension is omitted, i.e.\ $\h{\gamma} \equiv 0$. 
The third example of a liquid bridge explores the role of isotropic surface tension in a membrane surrounding a thin-walled  cylinder  with vanishing material properties. 
This example is a good benchmark for membrane problems as an analytical solution is available. 
It is however not ideally suited to problems in surface elasticity as, by definition, we require the energetic surface to be the boundary of a volume.
The fourth example illustrates neo-Hookean type surface effects in a nanoscale plate with a realistic rough surface. 
This example is used to assess the performance of the implementation for relatively complex and realistic geometries. 

The constitutive relations are given in Summary~\ref{constiutive_relations}. 
Trilinear and bilinear elements are used in the volume and on the surface, respectively. 
This can be modified in the input parameter file.
These results have been discussed, and additional context provided, in  the companion paper \citep{Javili-in-review}. 
The objective here is to provide information on the performance and features of the numerical scheme.

\subsection{Cook's membrane: neo-Hookean boundary potential}

\begin{figure}[!ht]
\centering
\includegraphics[width = \textwidth]{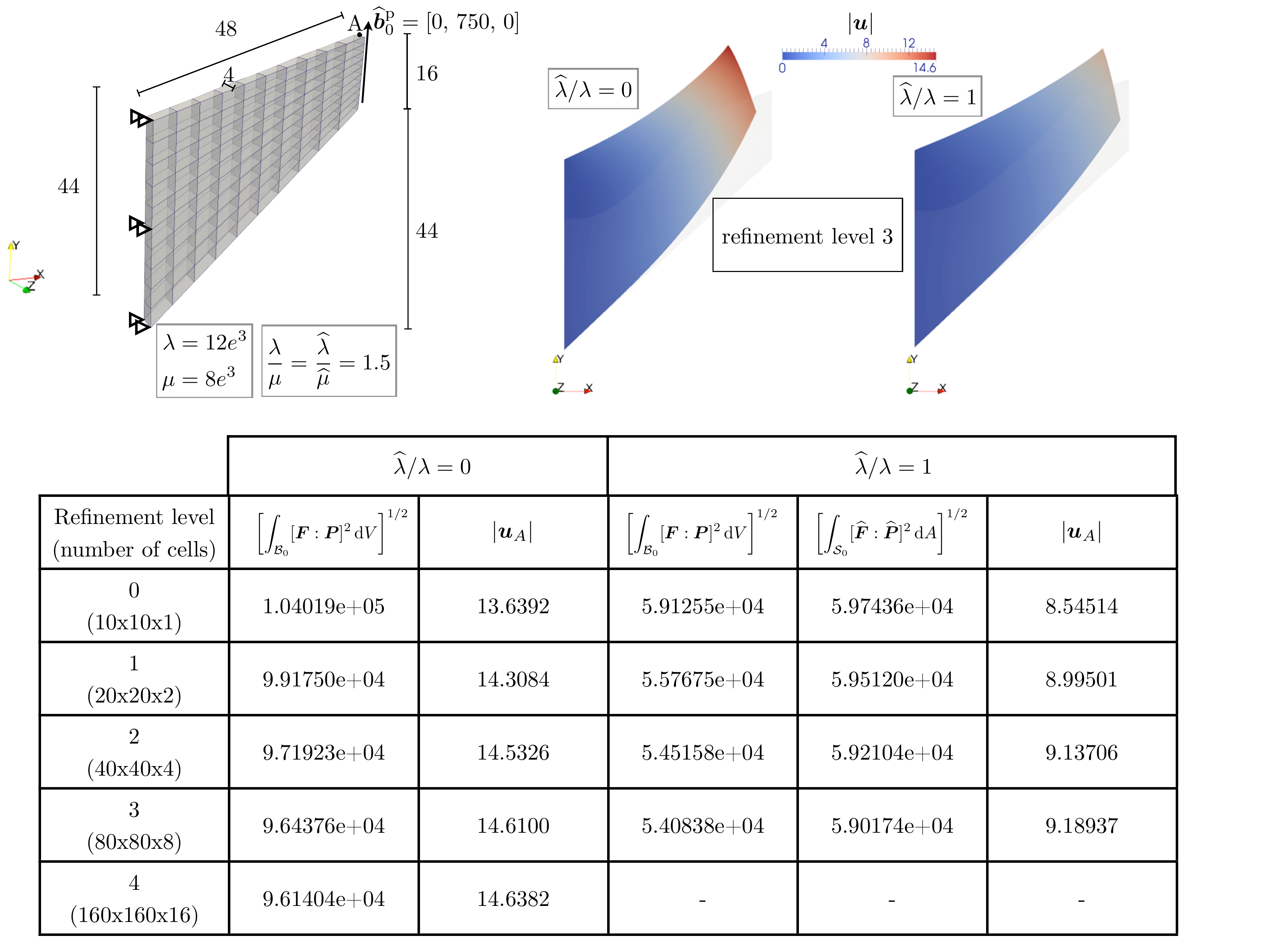}
 \caption{The initial 10x10x1 triangulation of the material configuration for the Cook's membrane. 
 The fixed material properties in the volume are given. 
 The final spatial configuration for the ratios $\h{\lambda} / \lambda = 0$ and $\h{\lambda} / \lambda = 1$ are shown. 
 The measured norms in the volume and on the energetic surface, and the displacement of point A, at four levels are refinement are recorded.
 Note that for $\h{\lambda} / \lambda = 1$ and a refinement level of 4, the simulation aborted. }
 \label{cooks}
\end{figure}

The Cook's membrane is a widely-used benchmark problem in solid mechanics. 
Consider the cantilever beam shown in Fig.~\ref{cooks}. 
The left face is fully fixed and a uniform traction applied to the right face. 
All of the remaining faces possess a surface energy (i.e.\ the back and front, and top and bottom faces).

The coarsest triangulation of the volume is obtained by uniformly dividing the horizontal and vertical edges into ten segments and creating a (non-affine) triangulation consisting of 10x10x1 cells. 
Subsequent triangulations are obtained by uniformly refining this initial one. 
The Neumann traction is applied in ten equal step. 

The material properties in the volume are given in Fig.~\ref{cooks}. 
The material properties of the energetic surface are varied while maintaining the ratio $\lambda / \mu = \h{\lambda} / \h{\mu} =1.5$.  
The following two norms in the volume and on the surface
\begin{align*}
\lrb{\int_{\c{B}_0} [\b{F} : \b{P}]^2  \, \d V}^{1/2}
&& \text{and} &&
\lrb{\int_{\c{S}_0} [\h{\b{F}} : \h{\b{P}}]^2  \, \d A}^{1/2} \, ,
\end{align*}
and the magnitude of the displacement of the point midway between the front and back faces on the edge common to the right and top faces, labelled $A$, for various levels of mesh refinement is also shown in Fig.~\ref{cooks}. 

The presence of the energetic surface significantly reduces the amount of deformation. 
A substantial amount of energy goes into deforming the energetic surface that would otherwise go into deforming the volume.
Numerical problems are encountered in the simulation on the finest mesh (160x160x16) when the energetic surface is present. 
The simulation aborts in the final time step due to a negative Jacobian determinant at a quadrature point. 
The reason for this appears to be due to the inherent length scale that the finite element implementation of surface elasticity theory introduces. 
The surface area to volume ratio of a cell (assumed to be a cube) scales as $4/h$ where $h$ is the edge length. 
If all the surfaces of the cell are assumed energetic, the contribution from the energetic surfaces of a volume cell will 64 times greater (relative to the volume contributions)  for the finest mesh relative to the coarsest. 
This can clearly introduce numerical problems and is similar to the issues that arise in heterogeneous material with vastly different material properties. 
This apparent deficiency is currently under investigation. 
In related work, the influence of the spatial discretization  in conjunction with the use of negative surface parameters, which have been reported in the literature, is assessed in \citet{Javili2012c}.

Although the results are clearly converging upon mesh refinement the presence of an extremely high stress concentration in the cells along the upper left edge of the domain limits the convergence rate.

\subsection{Nanowire: neo-Hookean boundary potential}

\begin{figure}[!ht]
\centering
\includegraphics[width = \textwidth]{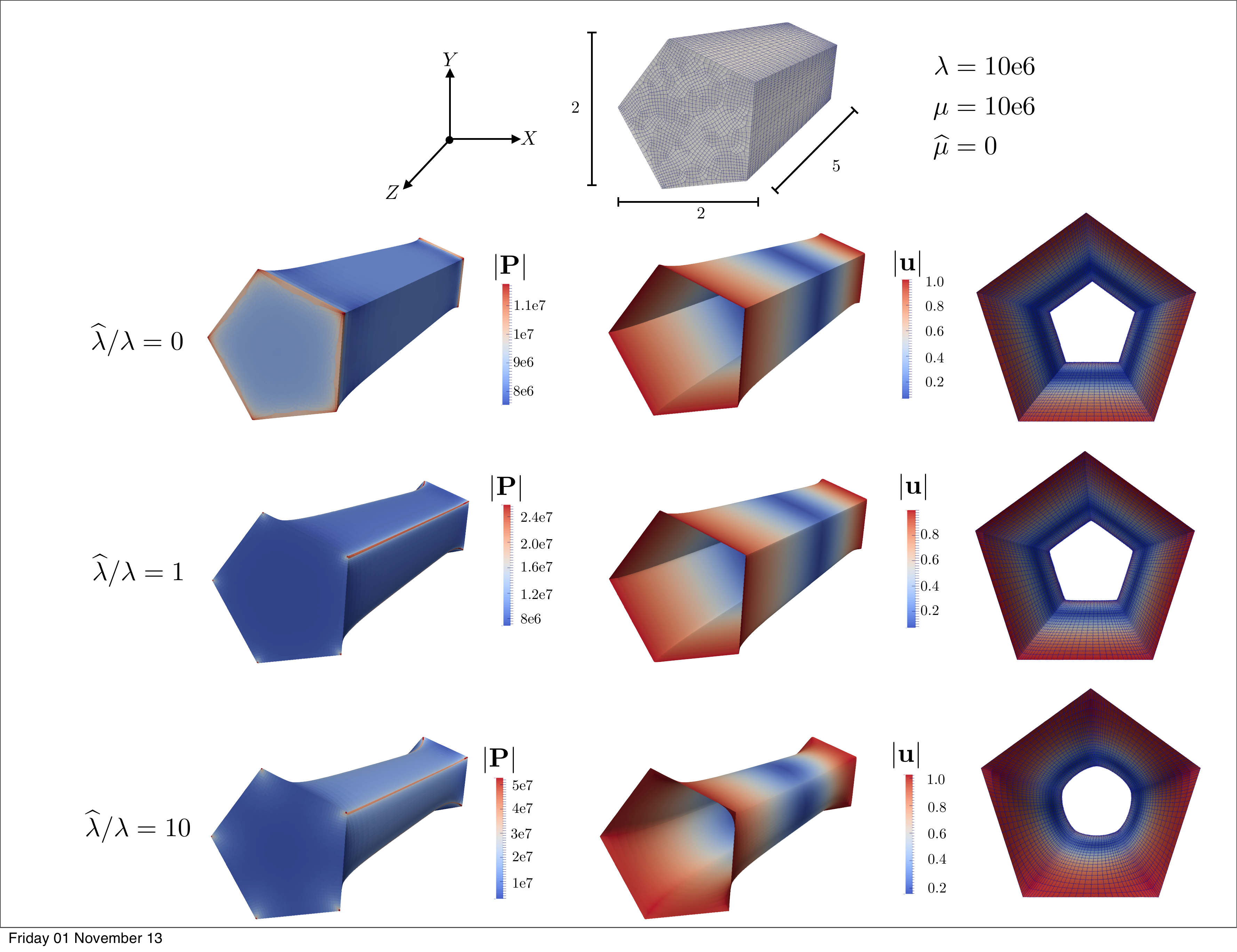}
 \caption{The triangulation of the material configuration for the nanowire. 
 The fixed material properties are given. 
 The final deformed (spatial) configuration of the volume and the surface for three different ratios of $\widehat{\lambda} / \lambda$ are shown.}
 \label{wire}
\end{figure}

One application of surface elasticity theory is to  describe surface effects in nanowires \citep[see e.g.][]{He2008, Yun2009, Yvonnet2011}. 
Consider the benchmark example shown in Fig.~\ref{wire}. 
The front and back pentagonal faces of the wire are prevented from displacing in the $X$ and $Y$ directions. 
The wire is extended in the $Z$-direction by an amount 2 (i.e.\ 40\% of the original length). 
The unconstrained surfaces on the side of the wire possess a surface energy. 

The triangulation of the volume, provided in {\tt nanowire\_fine.inp}, is more refined at the intersection of the faces on the surface and towards the front and back faces. 
This is done to better resolve the expected stress concentrations.
The volume and surface are discretized into 45 570 and 4500 elements, respectively. 
The prescribed deformation is applied uniformly in 10 steps.
An alternative strategy to generate a refined mesh would be adaptive mesh refinement.
A coarse mesh is also provided in the file {\tt nanowire\_coarse.inp}.

The Lam\'{e} parameters in the volume are fixed at the values given in Fig.~\ref{wire}. 
Similarly, the neo-Hookean energetic surface is characterised by the surface Lam\'{e} parameters 
$\widehat{\lambda}$ and $\widehat{\mu}$.
The surface shear modulus $\widehat{\mu}$ is set to zero and the ratio $\widehat{\lambda} / \lambda$ varied. 

The response in the absence of a surface energy is given by choosing $\widehat{\lambda} / \lambda = 0$ and is shown in  Fig.~\ref{wire}. 
The stress in the volume concentrates at the corners on the front and back faces. 
The initially pentagonal cross section reduces uniformly in size along the length. 

A surface energy is then assigned to the external surface. 
The stress in the volume $\b{P}$ concentrates along the lines forming the intersections of the external surfaces as  the value of $\widehat{\lambda}$ increases. 
Increasing the surface energy causes the resulting deformed cross section to tend to circular, thus increasing the stress in the volume in regions where the surface curvature is not smooth. 

The convergence history of the Newton scheme during the first timestep for $\h{\lambda} / \lambda = 1$ is given in Listing~\ref{convergence}. 
The various columns under the heading \verb|SOLVER STEP| correspond to, for each iteration, the application of the constraints (\verb|CST|), assembling the system matrix contributions from the volume and the surface (\verb|ASS_v| and \verb|ASS_s|, respectively), the solution of the linear system (\verb|SLV|), and updating the data stored at the continuum points in the volume and on the surface (\verb|UCP_v| and \verb|UCP_s|, respectively). 
The two columns labelled \verb|LIN_IT| and \verb|LIN_RES| record the number of iterations required by the linear solver and the final converged residual value, respectively. 
The final two columns give the norms of the global residual vector $\overset{\square}{\mathsf{R}}$ and the increment in the solution $\Delta \mathsf{u}$, defined in (\ref{newton_step}). 
The computation of these norms is restricted to unconstrained degrees of freedom. 
The Newton procedure continues until these two measures of the error decrease below the user-specified tolerance. 

\begin{lstlisting}[float,caption=The convergence history of the Newton scheme during the first timestep for the nanowire problem., float=h, label=convergence, style=Bash]
Timestep 1 @ 0.1s
_____________________________________________________________________________________
               SOLVER STEP                 |  LIN_IT  LIN_RES   | |R_NORM|  |dU_NORM|
_____________________________________________________________________________________
  0  CST  ASS_v  ASS_s  SLV  UCP_v  UCP_s  |     274  1.064e+02  1.000e+00  1.000e+00
  1  CST  ASS_v  ASS_s  SLV  UCP_v  UCP_s  |     795  2.873e+00  2.715e-02  3.603e-03
  2  CST  ASS_v  ASS_s  SLV  UCP_v  UCP_s  |     864  4.630e-02  4.425e-04  4.877e-05
  3  CST  ASS_v  ASS_s  SLV  UCP_v  UCP_s  |    1026  6.109e-04  6.134e-06  2.349e-07
  4  CST  ASS_v  ASS_s  SLV  UCP_v  UCP_s  |    1440  3.414e-08  3.187e-10  1.660e-10
 CONVERGED! 
_____________________________________________________________________________________
Relative errors:
	Solution: |dU_NORM|	1.660e-10
	Residual: |R_NORM|	3.187e-10
\end{lstlisting}

\subsection{Liquid bridge: isotropic surface tension effects} \label{sec:liquid_bridge}

\begin{figure}[!ht]
\centering
\includegraphics[width = \textwidth]{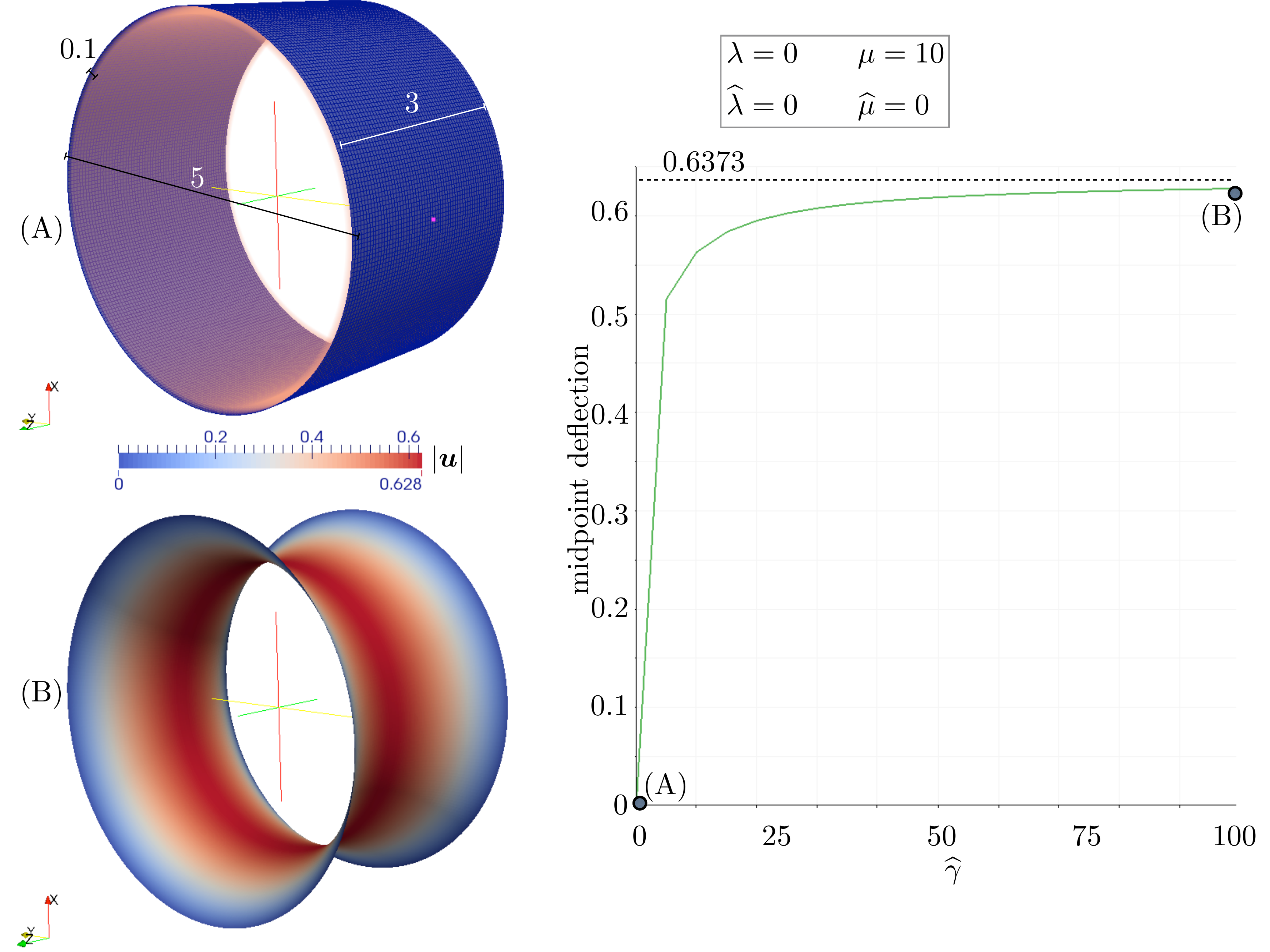}
 \caption{The liquid bridge in the material configuration (A) (volume and surface) and the final deformed state of the surface (B).
 The deflection of the midpoint of the surface versus the value of $\widehat{\gamma}$ is shown.}
 \label{surface_tension}
\end{figure}

A liquid bridge is a thin film suspended between two rigid side walls, as depicted in Fig.~\ref{surface_tension}. 
Surface tension acts in the material (initial) configuration to deform the surface so as to minimise its surface area. 
In order to model the liquid bridge using the approach to surface elasticity adopted here, a volume must be present. 
To that end, a neo-Hookean thin-walled cylinder (wall thickness = 0.1)  with Lam\'{e} parameters of $\lambda = 0$ and $\mu = 10$ is enclosed by the energetic surface. 
Ideally, the volume would have vanishing strength but numerically this results in a poorly conditioned system matrix.
The volume and surface triangulation contains 36 226 and 18 290 elements, respectively. 

The surface free energy, i.e.\ $\widehat{\gamma}$, is linearly increased over 20 equal steps to a value of $100$ and the displacement of a point on the middle of the surface monitored. 
The deflection of a midpoint on the surface for varying $\widehat{\gamma}$ is shown in Fig.~\ref{surface_tension}. 
The majority of the deformation occurs for $\widehat{\gamma} < 10$.
Thereafter, as $\widehat{\gamma}$ increases so the midpoint deflections converges to the analytical solution of 0.6373 \citep[see][]{Javili2010}. 

The numerical solution of this problem is somewhat sensitive. 
The rate at which the solution converges to the exact solution with increasing surface tension $\widehat{\gamma}$ is governed by the material properties of the volume, i.e.\ $\h{\mu}$, and the spatial discretization. 
The correct way to solve this problem robustly would be to disregard the volume and consider the membrane alone. 
This could be done relatively easily but would require the structure of the code to be modified significantly.

\subsection{Bending of a nanoscale plate with a rough surface}

\begin{figure}[!ht]
\centering
\includegraphics[width = \textwidth]{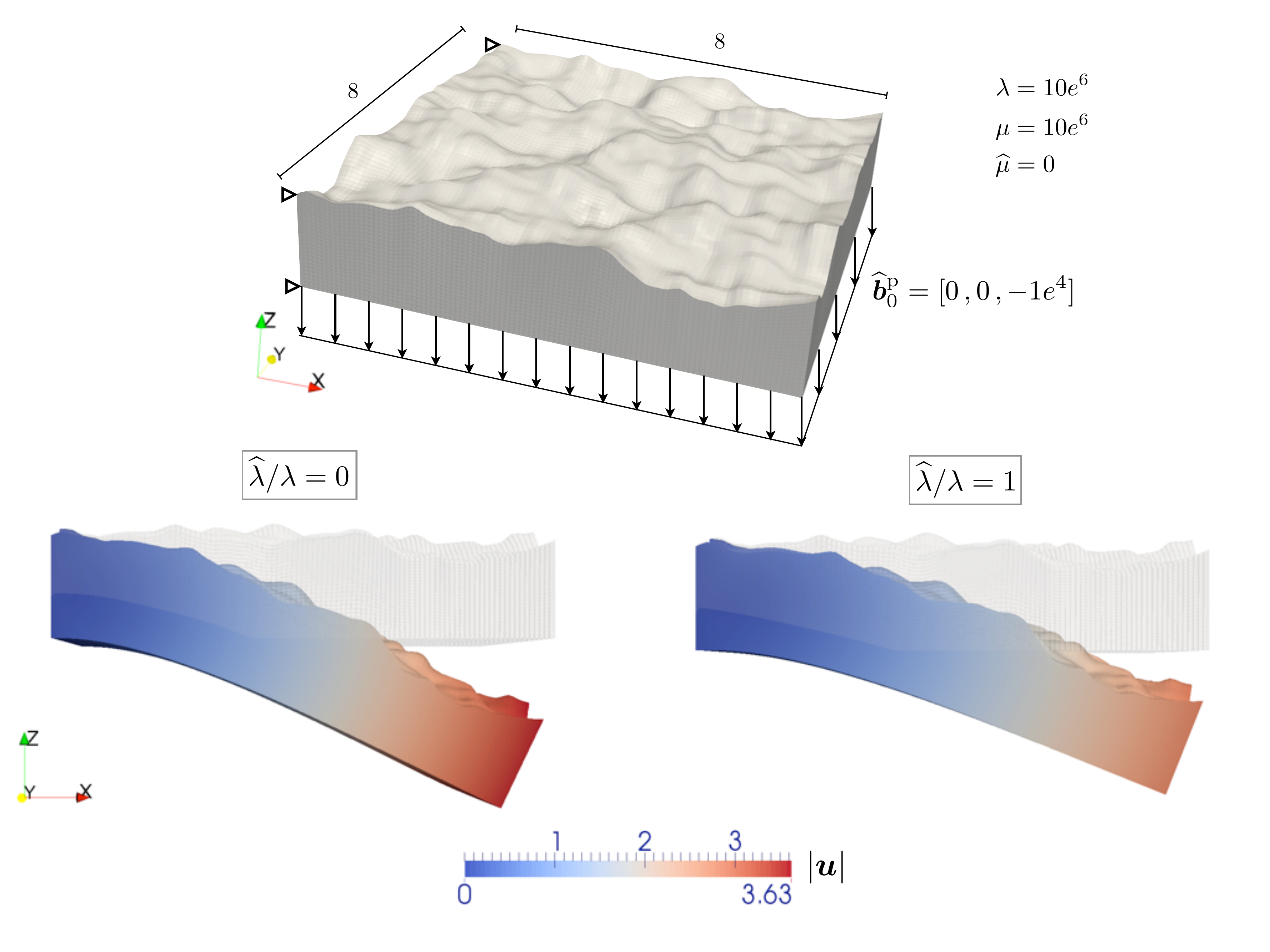}
 \caption{The triangulation and material configuration for the cantilever plate. 
 The fixed material properties are given. 
 The final spatial configuration of the bulk for two different ratios of $\widehat{\lambda} / \lambda$ are shown. 
 The magnitude of the displacement $\b{u}$ is plotted.}
 \label{plate}
\end{figure}

Surface roughness can have a significant, and often complex, influence on the response of nanoscale objects. 
Consider the cantilever plate shown in Fig.~\ref{plate}. 
The profile of the upper surface was produced using an open-source rough surface generation tool \citep{mysimlabs}.\footnote{The routine used, {\tt rsgene2D}, produces a Gaussian height distribution with an exponential auto-covariance. The input parameters were 100 divisions, a surface length of 2, a root mean square height of 0.05, and an (isotropic) correlation length of $0.25$. The surface was then scaled uniformly to have a surface length of 8.}
The rough upper surface is assumed energetic. 
All other surfaces are planar and standard. 
The Lam\'{e} parameters in the bulk are fixed at the values given in Fig.~\ref{plate}. 
The surface shear modulus $\widehat{\mu}$ is set to zero and the ratio $\widehat{\lambda} / \lambda$ varied. 
The left edge of the plate is fully fixed and a prescribed surface traction of $\h{\b{b}}\prs_0 = [0\, , 0\, , -1e^4]$ is imposed incrementally (in 10 uniform steps) on the lower face. 
The bulk and surface triangulation contains 250 000 and 10 000 elements, respectively (a coarser triangulation is provided in the accompanying online documentation).

The response of the plate to the loading for $\widehat{\lambda} / \lambda = 0$ and $\widehat{\lambda} / \lambda = 1$ is compared in Fig.~\ref{plate}. 
The plate with the energetic rough surface deforms considerably less. 

The example of the cantilever plate with the rough surface provides a good test of the robustness and efficiency of the numerical implementation. 
The Newton scheme exhibited the quadratic convergence associated with a consistently derived tangent. 
A summary of the time spent in the various sections of the code is given in Listing~\ref{timing}. 
It is clear that for problems of a reasonable size, the bottleneck is the solution of the linear system which accounts here for 77.7\% of the simulation time.  
A distributed parallel implementation utilising an AMG preconditioner, in the spirit of the approach used by \citet{Frohne2013}, would be one way to address this issue.

\begin{lstlisting}[float,caption=Summary of the time spent in the various parts of the code for the cantilever plate problem. The simulation was performed on a dual processor 8-core Intel Xeon 2.4GHz CPU with 64GB RAM, float=h, label=timing, style=Bash]
+---------------------------------------------+------------+------------+
| Total wallclock time elapsed since start    | 6.375e+03s |            |
|                                             |            |            |
| Section                         | no. calls |  wall time | % of total |
+---------------------------------+-----------+------------+------------+
| Assemble system volume          |        55 | 2.667e+02s |  4.18e+00% |
| Assemble tangent surface        |        55 | 1.795e+00s |  2.82e-02% |
| Construct grid                  |         1 | 7.110e+00s |  1.12e-01% |
| Linear solver                   |        55 | 4.950e+03s |  7.76e+01% |
| Postprocess results             |        11 | 8.970e+02s |  1.41e+01% |
| Setup system                    |         1 | 8.555e+00s |  1.34e-01% |
| Update CP data                  |       110 | 2.083e+02s |  3.27e+00% |
+---------------------------------+-----------+------------+------------+
\end{lstlisting}

\clearpage

\section{Conclusion}

An efficient finite element scheme for the solution of problems in nonlinear surface elasticity was presented. 
The effectiveness of the scheme was demonstrated using three example problems.
The complete, documented source code has been provided.
The code also provides a template for problems in classical finite elasticity.

Various extensions of the approach are planned. 
The first is to perform adaptive mesh refinement using the tools provided by {\tt deal.II}.
The key challenge here is the projection of the data conveniently located at the level of the quadrature point between successive meshes. 
This should not be too problematic as, unlike plasticity, none of the governing (evolution) equations are defined pointwise. 

The second extension would be to improve the efficiency of the code by implementing distributed parallelism. 
The primary challenge would be to ensure that the independent distributed volume and surface meshes can communicate, via a map of the common degrees of freedom.
The use of better-suited parallel solvers and preconditioners is also under investigation. 
Based on the results presented in \citet{Frohne2013}, for similar problems in elastoplasticity, the choice of the AMG preconditioner and BiCGStab solver should provide a good starting point.

The tangent is recomputed during each iteration of the Newton scheme. 
This operation is not that expensive, relative to the solution of the linear system, but may not be necessary. 
Finally, the Newton scheme should be extended to include a line search algorithm and adaptive time stepping introduced. 

A way to handle the numerical problems that arise as the (energetic) surface area to volume ratio increases upon mesh refinement is currently under investigation.

\section*{Acknowledgements}

The developers of {\tt deal.II} and the broader community of users are gratefully thanked.
The support for this work provided by the National Research Foundation through the South African
Research Chair for Computational Mechanics, and the Cluster of Excellence EAM (Multiscale Modelling and Simulation),   is gratefully acknowledged.

\begin{appendices}
\section{Library information and source code}

The complete, documented source code can be found online at \url{www.cerecam.uct.ac.za/code/surface_energy/doc/html/index.html}.
The code was compiled with version 8.0.0 of the {\tt deal.II} library using the precompiled Mac OS X package, and on a Linux machine using version 8.1.pre.

The steps to follow to run the numerical examples are described in the online documentation. 
The procedure to compile the code is identical to the {\tt deal.II} tutorials. 

\section{Features of the implementation in {\tt deal.II}}

The structure of the implementation follows the majority of the \verb|deal.II| tutorial examples.
The following section explains the key namespaces, classes and methods in the implementation. 

\subsection{Namespaces}

\subsection*{{\tt namespace Surface\_Elasticity}}
The complete implementation is wrapped within this namespace.

\subsection*{{\tt namespace Parameters}}
The parsing and reading of the input parameters from the various parameter files are handled in this namespace.

\subsection*{{\tt namespace AdditionalTools}}
This namespace defines several operations from tensor algebra that are not part of the {\tt deal.II} library.

\subsection{Classes and methods}

\subsection*{{\tt class Time}}
Simple time management class used to incrementally advance the quasi-static problem.

\subsection*{{\tt class ContinuumPoint}}
As discussed in Section~\ref{sec:continuum_points}, this class is responsible for the handling the constitutive response. 
This class would need to be modified if a different type of material model was adopted. 

\subsection*{{\tt class Solid}}
This class contains the core routines to manage to the code. 
A {\tt Solid} encapsulates the volume and the surface and all the routines required to control the problem. 
The key methods within the class are:
\subsubsection*{{\tt run ()}} This is the method responsible for the high level control. 
It calls the routines to create the grid, solve the problem using the Newton scheme, and to output the results.

\subsubsection*{{\tt system\_setup ()}}
This method is used to initialise the system matrix and various other vectors, including the incremental solution and the right-hand side. 

\subsubsection*{{\tt assemble\_system\_volume () }}
This method is responsible for assembling the contributions from the volume to the system matrix and the right-hand side vector. 

\subsubsection*{{\tt assemble\_system\_surface () }}
This method is responsible for assembling the contributions from the surface to the system matrix and the right-hand side vector.

\subsubsection*{{\tt make\_constraints (...)}}
The Dirichlet constraints are assembled in this routine. 

\subsubsection*{{\tt solve\_nonlinear\_timestep (...)}}
This routine controls the Newton--Raphson procedure used to solve the problem. 
It continues the iterative procedure until a converged solution is obtained.

\subsubsection*{{\tt solve\_linear\_system}}
This routines solves the linear system of equations using the choice of solver and preconditioner specified in the parameter file.

\subsubsection*{{\tt output\_results ()}}
This routine is responsible for outputting the results in \verb|vtu| format for subsequent post-processing.
The user can choose to write the data associated with the quadrature points. 
Note, this can be time consuming for large problems.
The solution in the volume and on the surface are outputted in separate files. 

\subsubsection*{{\tt update\_volume\_cp\_incremental}}
Given an updated solution, this routine is responsible for updating the data stored at the quadrature points in the volume.

\subsubsection*{{\tt update\_surface\_cp\_incremental}}
Given an updated solution, this routine is responsible for updating the data stored at the quadrature points on the surface.

\end{appendices}


\bibliographystyle{plainnat}
\bibliography{Elibrary}

\end{document}